\documentclass[11pt]{amsart}

\usepackage{xcolor}
\usepackage[verbose]{hyperref}
\hypersetup{colorlinks=false,allbordercolors=blue,pdfborderstyle={/S/U/W 1}}

\usepackage{amssymb}
\usepackage{amsmath}
\usepackage{amsthm}

\theoremstyle{plain} 
    \newtheorem{theorem}{Theorem}[section]   
    \newtheorem{lemma}[theorem]{Lemma}
    \newtheorem{corollary}[theorem]{Corollary}

\theoremstyle{definition} 
    \newtheorem{definition}[theorem]{Definition}
    \newtheorem{example}[theorem]{Example}

\usepackage{tikz}
\usetikzlibrary{graphs}
          
\newcommand{\id}{e}                                      
\newcommand{\csg}[1]{\langle #1\rangle}       
\newcommand{\gen}[1]{\mathrm{gen}({#1})}   
\newcommand{\ord}[1]{\left\vert#1\right\vert}  
\newcommand{\oo}[1]{{\mathrm{o}(#1)}}         
\newcommand{\aut}[1]{\mathrm{Aut}(#1)}       

\newcommand{\cgp}[2]{{#1}^{ #2 }}  

\newcommand{\dpk}[2]{\vec{\mathcal{P}}_{#1}({#2})}
\newcommand{\pk}[2]{\mathcal{P}_{#1}({#2})}
\newcommand{\dqk}[2]{\vec{\mathcal{Q}}_{#1}({#2})}
\newcommand{\qk}[2]{\mathcal{Q}_{#1}({#2})}

\newcommand{\exclude}[1]{\widetilde{#1}}

\begin{document}

\markboth{Brian Curtin}
{Excluded Power Graphs}

%
%

\title{Excluded power graphs of groups}

\author{Brian Curtin}

\address{
Department of Mathematics and Statistics\\
University of South Florida\\
4202 E. Fowler Ave., CMC 342\\
Tampa, FL 33620, USA
}

\email{bcurtin@usf.edu}{}

\bigskip

\begin{abstract}
Let $\mathcal{X}$ be a set of integers greater than one. 
The $\mathcal{X}$-excluded power graph of a group $G$ has vertex set $G$ and an edge from $g$ to each power of $g$ other than itself provided that the power is not divisible by any element of $\mathcal{X}$.  
When $G=H\times K$ for groups $H$ and $K$ with coprime orders, excluding the prime factors of $|H|$ yields a power graph with a quotient consisting of multiple copies of a quotient of the power graph (no exclusions) of $K$.   Partial results for the semidirect product under the same conditions are given.
We describe groups whose $\mathcal{X}$-excluded power graphs consist of disjoint directed cliques.
\end{abstract}



\maketitle

\section{Introduction}

There is a long and active tradition of using graphs to study groups.   One such graph is the power graph of a group.  Power graphs were introduced for semigroups \cite{MR1900273} and groups \cite{MR1777663} by Kelarev and Quinn.  Undirected power graphs were introduced by Chakrabarty et al.~\cite{MR2511776}.   A number of graphs constructed from groups, including power graphs of several flavors, are surveyed  in \cite{MR4346241}, and a hierarchy of a number of such graphs is studied in \cite{MR4798847,MR4426775}.

\begin{definition}
\label{def:powergraph}
Let $G$ be a group.  
The {\em directed power graph} $\dpk{}{G}$ and  
the {\em(undirected) power graph} $\pk{}{G}$ of $G$ 
each has vertex set $G$ and, respectively, have 
a directed edge from $g$ to $h$  whenever 
       $h\not=g$ and $h=g^k$ for some positive integer $k$ and 
an undirected edge between distinct $g$ and $h$ whenever  
      $h\not=g$ and either $h=g^k$ or $g=h^k$ for some positive integer $k$.
\end{definition}

Power graphs are related to the cyclic subgroups of a group. 
Indeed, in the directed power graph there is an edge from $g$ to $h$ when the cyclic subgroup generated by $h$ is contained in the cyclic subgroup generated by $g$,  and 
in the undirected power graph  two elements are adjacent when the cyclic subgroup generated by one is contained in the cyclic subgroup generated by the other.

Several variations of power graphs have been studied. Some remove edges that arise from basic group properties rather than the properties  of a specific group.
The {\em proper power graph} deletes the identity and its incident edges \cite{MR3366577}.
The {\em reduced power graph} omits edges of the power graph between group elements that generate the same cyclic subgroup (corresponding to strict containment of the cyclic subgroups that they generate) \cite{MR4269515}.
Some variations of power graphs focus on other interactions of cyclic subgroups.
The {\em enhanced power graph} connects distinct group elements which are in a common cyclic subgroup \cite{MR4446436}.
The {\em intersection power graph} joins group elements when the cyclic subgroups they generate intersect nontrivially \cite{MR3790706}.
The {\em trivial intersection power graph} joins group elements when the cyclic subgroups they generate intersect trivially \cite{MR4849790}.
The {\em enhanced power coprime graph} joins group elements with coprime order when they belong to a common proper cyclic subgroup \cite{MR4958256}.
Major themes in the literature include 
   describing graph theoretic and spectral properties of various power graphs arising from groups with some property, 
   characterizing the groups whose power graphs satisfy various graph theoretic conditions, and 
   understanding the relationship between power graphs and other graphs defined on groups.  
Power graphs of groups are surveyed in \cite{MR3145411,MR4310374,MR4446436}.  

In this paper we introduce another variation of the power graph of a group.   We refine the information provided by the power graph of a group and to isolate certain of its edges by excluding edges that arise from certain powers of elements.

\begin{definition}
\label{def:excludedpowergraph}
Let $G$ be a group.  Let $\mathcal{X}$ be a subset of integers greater than one.   
The {\em directed $\mathcal{X}$-excluded power graph}  $\dpk{\exclude{\mathcal{X}}}{G}$ and 
the {\em (undirected) $\mathcal{X}$-excluded power graph}  $\pk{\exclude{\mathcal{X}}}{G}$ of $G$
each has vertex set $G$ and, respectively, have an edge from $g$ to a distinct $h$ whenever $h=g^k$ for some positive integer $k$ not divisible by any element of $\mathcal{X}$ 
and an undirected edge between distinct $g$ and $h$ whenever either $h=g^k$ or $g=h^k$ for some positive integer $k$ not divisible by any element of $\mathcal{X}$. 
\end{definition}

The edges of an excluded power graph are among these of the power graph.  We use excluded power graphs to help describe some substructures of power graphs of certain direct and semidirect products of groups which are inherited from the factors.   
In \cite{MR3754806}, Bhuniya and Mukherjee introduced a generalized graph product with the property that the power graph of the direct product of groups is isomorphic to the generalized product of their power graphs.   
Other papers which examine power graphs of direct products include \cite{MR4764545,MR4479798,MR4525231}.
Describing power graphs of semidirect products of groups is more challenging, so we welcome the small insight that excluded power graphs offer when the factors have coprime order.

In Section \ref{sec:basics}, we describe a few basic properties of $\mathcal{X}$-excluded power graphs, although our main focus is the case where $\mathcal{X}$ is a set of primes.   
In Section \ref{sec:hallexcluded},  we consider  $\mathcal{X}$-excluded power graphs of semidirect products of $H\rtimes K$, where $H$ and $K$ have coprime orders and $\mathcal{X}$ consists of the prime factors of $|H|$.  
In particular, we show that when $G=H\times K$ for groups $H$ and $K$ with coprime orders, excluding the prime factors of $|H|$ has quotient power graph consisting of multiple disjoint copies of the quotient of the power graph of $K$.  
Finally, in Section \ref{sec:disjointcliques} we describe those groups whose directed and undirected $\pi$-excluded power graphs are a disjoint union of cliques for some set of primes $\pi$: 
In the directed case, the group must be a $\pi$-group, and in the undirected case, the group must be the semidirect product of a cyclic group of prime power order by a $\pi$-group.

\section{Background}\label{sec:graphs}

We first recall some graph terminology.
A {\em directed graph} or {\em digraph} consists of a set of {\em vertices} and a set of  {\em (directed) edges} comprised of ordered pairs of vertices.  A directed edge $(v, w)$ is  {\em from $v$ to $w$}, 
$v$ is an {\em in-neighbor} of a vertex $w$, and $w$ is an {\em out-neighbor} of $v$. 
A {\em loop} is an edge from a vertex to itself.           
A subset of vertices of a digraph form a {\em directed clique} whenever there is a directed edge from each vertex of the clique to every other vertex of the clique.    
A {\em (simple) graph} consists of  a set  of vertices and a set of {\em edges} comprised of unordered pairs of distinct vertices.  An edge $\{v, w\}$ is  {\em between $v$ and $w$}, and  $v$ and $w$ are {\em neighbors} of one another.

Let $\Gamma_1=(V_1, E_1)$ and $\Gamma_2=(V_2, E_2)$ be (di)graphs.  A {\em graph isomorphism} from $\Gamma_1$ to $\Gamma_2$ is a bijection  $\psi:V_1\rightarrow V_2$ such that $(v,w)\in E_1$ if and only if $(\psi(v), \psi(w))\in E_2$. 
We say that $\Gamma_1$ and $\Gamma_2$ are {\em isomorphic} whenever there is a graph isomorphism from one to the other, and in this case we write $\Gamma_1\cong \Gamma_2$.

 A graph $G$ is {\em connected} if there exists a path between every pair of vertices. A {\em (connected) component} of $G$ is a subgraph that is connected and is not properly contained in any other connected subgraph of $g$.
A {\em weakly connected component} of a digraph is a connected component in the undirected graph formed by replacing each directed edge with an undirected edge. 
%
An {\em induced subgraph} of a (di)graph $\Gamma$ is a graph whose vertices are a subset of those of $\Gamma$  and whose edges are all those of $\Gamma$ which join two vertices in the subset.  
The power graph (of any flavor) of a subgroup is the induced subgraph of the power graph (of the same flavor) of the group on that subgroup. 

The {\em quotient graph} of a (di)graph relative to a partition of its vertices takes the cells of the partition as its vertices with an edge from one cell to another whenever there is an edge from a vertex of the first cell to a vertex of the second cell.  
A partition $\{V_1, V_2, \ldots, V_k\}$ of the vertices of a (di)graph is called an {\em equitable partition} if for every pair of indices $i$ and $j$, every vertex in $V_i$ has the same  number of (in- and out-)neighbors in $V_j$.

\begin{definition}
Let $G$ be a group, and let $\mathcal{C}$ be the set of cyclic subgroups of $G$. 
The {\em directed quotient power graph} $\dqk{}{G}$ and 
the {\em (undirected) quotient power graph} $\qk{}{G}$  of $G$ 
each has vertex set $\mathcal{C}$ and, respectively, have 
a directed edge from each cyclic subgroup to each of its proper subgroups and 
an undirected edge joining two distinct cyclic subgroups whenever one is a proper subgroup of the other. 
Label the directed edge from $C$ to $D$ by the index of $D$ in $C$.
Let $\mathcal{X}$ be a set of integers greater than one.
The {\em directed $\mathcal{X}$-excluded quotient power graph} $\dqk{\exclude{\mathcal{X}}}{G}$ and 
the {\em (undirected) $\mathcal{X}$-excluded quotient power graph} $\qk{\exclude{\mathcal{X}}}{G}$ of $G$ 
are the subgraphs of $\dqk{}{G}$ formed by deleting the edges of the corresponding power graph with labels divisible by some element of $\mathcal{X}$.
\end{definition}

With a bit of elementary group theory we describe the quotient power graph as a quotient graph.
A {\em homomorphism} from one group $G_1$ to another group $G_2$ is a map $\phi:G_1\rightarrow G_2$ satisfying $\phi(gh)=\phi(g)\phi(h)$ for all $g$, $h\in G_1$.  
An {\em isomorphism} from one group $G_1$ to another group $G_2$ is a bijective homomorphism from $G_1$ to $G_2$.  We say that $G_1$ and $G_2$ are isomorphic whenever there is an isomorphism from one to the other, 
and in this case we write $G_1\cong G_2$.  An {\em automorphism} of a group $G$ is an isomorphism from $G$ to itself.  We write $\aut{G}$ for the automorphism group of $G$.

 Let $G$ be a group with identity $\id$. Let $\ord{G}$ denote the {\em order} (cardinality) of $G$.  
 For $g\in G$, let $\csg{g}$ denote the cyclic subgroup of $G$ generated by $g$.
  The {\em order $\oo{g}$ of $g$} is both $\ord{\csg{g}}$ and 
                   the least positive integer $t$ such that $g^{t}=\id$ (or $\infty$ if no such $t$ exists).  
%
Let $C$ be a cyclic group of finite order $n$.    Let $\gen{C}$ denote the set of generators of $C$.   
Given $g\in\gen{C}$,  $\gen{C}= \{ g^m : \gcd(m, n)=1, 1\leq m<\oo{g}\}$.
Observe that  $\gen{C}$ has size  $\varphi(n)$, where $\varphi$ is the Euler totient. 
For each nonnegative integer $k$,  define $\pi_k:h\mapsto h^k$ $(h\in C)$.
Then $\pi_k$ is an automorphism of $C$ if and only if 
         $\pi_k$ permutes generators $C$ if and only if 
         $k$ and $n$ are coprime.
In fact, $\aut{C} =\{\pi_k : \gcd(k, n)=1\}$ and 
            $\ord{\aut{C}}=\varphi(n)$. 

The quotient power graph is a graph quotient relative to the partition of group elements whose cells consist of the generators of cyclic subgroups, with the loops deleted. 

\begin{lemma}  (c.f.~\cite[Lemma 3.3]{MR3743245})
Let $G$ be a group, and let $\mathcal{C}$ be the set of cyclic subgroups of $G$. 
In  $\dpk{}{G}$,  $\gen{C}$ is a directed clique, and 
every element of $\gen{C}$ has the same in- and out-neighbors not in $\gen{C}$.
In particular,  $\{\gen{C} : C\in\mathcal{C}\}$ is an equitable partition of $\dpk{}{G}$ and $\pk{}{G}$.
Moreover, $\dqk{}{G}$ and $\qk{}{G}$ are, respectively, quotient graphs of  $\dpk{}{G}$ and $\pk{}{G}$ with respect to this equitable partition.
\end{lemma}

The basic properties of finite cyclic groups give the following.

\begin{lemma}
Let $G$ be a group, and let $\mathcal{C}$ be the set of cyclic subgroups of $G$. 
Pick $C\in\mathcal{C}$, and let  $n=\ord{C}$.   
For $k\in \mathbb{Z}$, define $\cgp{C}{k} = \{h^k:h\in C\}$.
\begin{enumerate}
\item For any positive integer $m$,  $\cgp{C}{m} = \cgp{C}{\gcd(m, n)}$.  
         In particular, $\cgp{C}{m} = C$ if and only if $\gcd(m,n)=1$.
\item For each divisor $d$ of $n$ other than 1, there is an edge from $C$ to $\cgp{C}{d}$ in $\dqk{}{G}$. 
Every edge of $\dqk{}{G}$ arises in this way.
The edge from $C$ to $\cgp{C}{d}$ is labeled $d=|C: \cgp{C}{d}|$.
\end{enumerate}
\end{lemma}

We will use one  elementary number theory fact.

\begin{lemma}\label{lem:exlucdeddivisors}
Fix a positive integer $n$ and a finite set $\mathcal{X}$ of integers greater than $1$.
For any integer $m$ $(1\leq m \leq n)$, the following are equivalent.
\begin{enumerate}
\item No element of $\mathcal{X}$ divides both $m$ and $n$.
\item There is a nonnegative integer $c$ such that  for all $x\in \mathcal{X}$, we have  $m+cn\not\equiv 0 \pmod x$.
\item The arithmetic sequence $m$, $m+n$, $m+2n$, $m+3n$, \ldots\, contains a number not divisible by any element of $\mathcal{X}$.
\end{enumerate}
\end{lemma}

\begin{proof} 
Suppose that for all nonnegative integers $c$, $m+cn\equiv 0\pmod x$ for all $x\in \mathcal{X}$.
Let $L$ be the least common multiple of the elements of $\mathcal{X}$. Then $L\mid m+ 0\cdot n$ and $L\mid m+1\cdot n$, so $L\mid m, n$.  In particular, every element of $\mathcal{X}$ divides both $m$ and $n$. 
Thus, (i) implies (ii). Conversely, if $x\in\mathcal{X}$ divides both $m$ and $n$, then $x$ divides $m+cn$ for all integers $c$.  Thus, (ii) implies (i).  
Note that (iii) is a restatement of (ii).
\end{proof} 

\section{Excluded Power Graphs}\label{sec:basics}

In this section we discuss edges of excluded power graphs.  
In the directed power graph of a group $G$, 
     the out-neighbors of $g\in G$ are $\{ g^k : 2\leq k \leq \oo{g}\}$. 
An exponent $k$ may be divisible by some element(s) of a set $\mathcal{X}$ of integers while $k+c\oo{g}$ is not, yielding an edge from $g$ to $g^k$ in $\dpk{\exclude{\mathcal{X}}}{G}$ despite appearing to be excluded based on the value of $k$.

\begin{lemma}\label{lem:weakexcludededge}
Let $G$ be a group, and pick $g\in G$ with finite order.   
Let $\mathcal{X}$ be a finite set of integers greater than 1, and let $k$ be an integer $2\leq k \leq \oo{g}$.
Then there is an edge from $g$ to $g^k$ in $\dpk{\exclude{\mathcal{X}}}{G}$ if and only if 
no element of $\mathcal{X}$ divides $|\csg{g}:\csg{g^k}|=\gcd(k, \oo{g})$.
\end{lemma}

\begin{proof} 
Fix $k$ $(2\leq k \leq \oo{g})$.  Take $m=k$ and $n=\oo{g}$ in Lemma \ref{lem:exlucdeddivisors} to see that there is some $\ell$ with $g^\ell=g^k$ and $\ell$ not divisible by any element of $\mathcal{X}$. 
 \end{proof} 

\begin{corollary}\label{cor:genarecliques}
Let $G$ be a group, and let $C$ be a finite cyclic subgroup of $G$. 
Let $\mathcal{X}$ be a finite set integers greater than 1.
Then $\gen{C}$ is a directed clique in $\dpk{\exclude{\mathcal{X}}}{G}$.
\end{corollary}

\begin{proof} 
Suppose $g\in\gen{C}$.  Then $\gen{C}=\{g^k:\gcd(k,\ord{C})=1\}$.  
When $\gcd(k,\ord{C})=1$, there is no element of $\mathcal{X}$ dividing both $k$ and $\ord{C}$,
so there is a directed edge from $g$ to $g^k$ by Lemma \ref{lem:weakexcludededge}.
\end{proof} 

\begin{lemma}
\label{lem:reversecontain}
Let $G$ be a group. Let $\mathcal{X}$ and $\mathcal{Y}$ be  finite sets of integers greater than 1.  
If $\mathcal{X}\subseteq \mathcal{Y}$, then the edge set of  $\dpk{\exclude{\mathcal{X}}}{G}$ contains  
the edge set of  $\dpk{\exclude{\mathcal{Y}}}{G}$.
\end{lemma}

\begin{proof}
If there is an edge from $g$ to $h$ in  $\dpk{\exclude{\mathcal{Y}}}{G}$, then $h=g^k$ for some $k$ not divisible by any $y\in Y$.
Now $h=g^k$ for some $k$ not divisible by any $x\in X\subseteq Y$, so there is an edge from $g$ to $h$ in $\dpk{\exclude{\mathcal{X}}}{G}$.
\end{proof}

\begin{corollary}
Let $G$ be a group. Let $\mathcal{X}$ and $\mathcal{Y}$ be  finite sets of integers greater than 1. 
 \begin{enumerate}
\item \label{lem:intersectionexcludededges}
         The edge set of $\dpk{\exclude{\mathcal{X}\cup\mathcal{Y}}}{G}$ is the intersection of
         the edge sets $\dpk{\exclude{\mathcal{X}}}{G}$ and $\dpk{\exclude{\mathcal{Y}}}{G}$. 
         
\item \label{lem:unionexcludededges}
         The edge set of $\dpk{\exclude{\mathcal{X}\cap \mathcal{Y}}}{G}$ contains the edge sets of $\dpk{\exclude{\mathcal{X}}}{G}$ and $\dpk{\exclude{\mathcal{Y}}}{G}$.

\end{enumerate}
\end{corollary}

\begin{proof}
By Lemma \ref{lem:reversecontain}, the edges of $\dpk{\exclude{\mathcal{X}\cup\mathcal{Y}}}{G}$ are edges of both $\dpk{\exclude{\mathcal{X}}}{G}$ and $\dpk{\exclude{\mathcal{Y}}}{G}$.
If there is an edge from $g$ to $h$ in both $\dpk{\exclude{\mathcal{X}}}{G}$ and $\dpk{\exclude{\mathcal{Y}}}{G}$,  then $h=g^k$ for some $k$ $(1<k\leq \oo{g})$. 
By Lemma \ref{lem:weakexcludededge}, no element of $\mathcal{X}$ and no element of $\mathcal{Y}$ divides $\gcd(k, \oo{g})$.  
Thus, no element of $\mathcal{X}\cup\mathcal{Y}$ divides $\gcd(k, \oo{g})$, so there is an edge from $g$ to $h=g^k$ in $\dpk{\exclude{\mathcal{X}\cup\mathcal{Y}}}{G}$.  
Thus, (i) holds.
Lemma \ref{lem:reversecontain} gives that the edges of $\dpk{\exclude{\mathcal{X}}}{G}$ and $\dpk{\exclude{\mathcal{Y}}}{G}$ are edges of $\dpk{\exclude{\mathcal{X}\cap \mathcal{Y}}}{H}$, 
so (ii) holds.
\end{proof}

\begin{example}
Equality need not hold in Lemma \ref{lem:unionexcludededges}.
Let $G$ be the cyclic group of order $12$, and $g$ be an element of order 12, and let $h=g^6$.  
There is no edge from $g$ to $h$ in either $\dpk{\exclude{\{2\}}}{G}$ or  $\dpk{\exclude{\{3\}}}{G}$. 
However, there is an edge from $g$ to $h$ in  $\dpk{\exclude{\{2\}\cap \{3\}}}{G} = \dpk{}{G}$.
\end{example}

We are interested in excluding prime divisors of the group order.  

\begin{lemma}
Let $G$ be a group, and let $H$ be a finite subgroup of $G$.
Let $\mathcal{X}$ and $\mathcal{Y}$ be  finite sets of prime numbers.
\begin{enumerate}
\item \label{lem:excludednondiv}
         Suppose no prime divisor of $\ord{H}$ is in $\mathcal{X}$. 
         Then $\dpk{\exclude{\mathcal{X}}}{H} = \dpk{}{H}$.
\item \label{lem:excludeallprimedivisors}
         Suppose every prime divisor of $\ord{H}$ is in $\mathcal{X}$.
         Then  $\dpk{\exclude{\mathcal{X}}}{H}$ is the disjoint union of directed cliques
         on the generators of each cyclic subgroup of $H$.
\end{enumerate}
\end{lemma}

\begin{proof} 
In (i), there is no element of $\mathcal{X}$ dividing $\oo{g}$ for any $g\in H$.
Thus, for any $k$ $(2\leq k \leq \oo{g})$ there is an edge from $g$ to $g^k$ in $\dpk{\exclude{\mathcal{X}}}{H}$
by Lemma \ref{lem:weakexcludededge}.

In (ii),  Corollary \ref{cor:genarecliques} implies that the generators of $C$ comprise a directed clique. Pick $g\in H\setminus\{\id\}$. Note that $\oo{g} \mid \ord{H}$, so it is a product of primes in $\mathcal{X}$.  
Lemma \ref{lem:weakexcludededge} gives there is no edge from $g$ to $g^k$ for any $k$ not coprime to $\oo{g}$. Thus, there are no edges other than those in the directed cliques, giving (ii).
\end{proof} 
  
\begin{example}\label{example:A5excluded}
The alternating group  $A_5$ consists of 
     the identity,  
     15 products of two disjoint two-cycles, 
     20 three-cycles (giving ten pairs of inverse elements which generate the same cyclic subgroup), and 
     24 five-cycles (giving six sets consisting of four five-cycles which generate the same cyclic subgroup).   
In $\dpk{\exclude{\{2\}}}{A_5}$, $\dpk{\exclude{\{3\}}}{A_5}$, and $\dpk{\exclude{\{5\}}}{A_5}$ we have edges among the five-cycles in each cyclic subgroup as
     $(abcde)^2 = (abcde)^7 = (acebd)$, 
     $(abcde)^3=(abcde)^8=(adbed)$, and 
     $(abcde)^4=(abcde)^9 = (aedcb)$. 
By Lemma \ref{lem:intersectionexcludededges}, the edges of $\dpk{\exclude{\{2,3\}}}{A_5}$ are those edges in both $\dpk{\exclude{\{2\}}}{A_5}$ and $\dpk{\exclude{\{3\}}}{A_5}$.  
\end{example}  

Although we focus on excluding only finite sets of integers (prime divisors of the group order), we comment on one family of infinite excluded sets. 
The  {\em iterated $e^\mathrm{th}$-power digraph} (also often referred to as a power graph) of a group has an edge from each group element $g$ to $g^e$. 
Iterated power graphs on $\mathbb{Z}_n$ have been extensively studied for their connection to number theory \cite{MR2167719,MR2165548,MR2059267} and applications, such as pseudorandom number generation \cite{zbMATH03972053,MR2165548}.   
Let $p$ be a prime and let $\mathcal{X}$ be all integers greater than 1 other than $p$.  Then $\dpk{\exclude{\mathcal{X}}}{G}$ is the iterated $p^\mathrm{th}$-power digraph.


\section{Semidirect Product of Groups of Coprime Order}\label{sec:hallexcluded}

Let $\pi$ be a set of primes.  A group $H$ is a {\em $\pi$-group} whenever all prime factors of $\ord{H}$ lie in $\pi$.  
A  subgroup $H$ of a group $G$ is a {\em Hall $\pi$-subgroup} whenever it is a $\pi$-group and the index of $H$ in $G$ is not divisible by any prime in $\pi$.  
In this section we comment on $\pi$-excluded power graphs in groups with a normal Hall $\pi$-subgroup.  The structure of such groups can be described as follows. 
A subgroup $H$ of a group $G$ is {\em complemented} whenever there exists a subgroup $K$ of $G$ such that $G=HK$ and $H\cap K=\{\id\}$. In this case, $K$ is called a {\em complement} of $H$. 
The Schur-Zassenhaus theorem asserts that every normal Hall subgroup $H$ has a complement $K$.
Moreover, $G= H\rtimes K$ when $H$ is normal and $K$ is a complement of $H$ in $G$ with $k\in K$ acting on $H$ by conjugation by $k$. 
See, for instance, Rotman \cite{MR1307623}.

\begin{theorem}\label{thm:complementedhall}
Let $\pi$ be a set of primes, and suppose $G$ is a finite group with a normal Hall $\pi$-subgroup $H$.
Let $K$ be a complement of $H$ in $G$, so $G=H\rtimes K$.   
Pick $a\in H$ and $b\in K$. Let $J=\{a'b': a'\in \csg{a}, b'\in \csg{b}\}$.
\begin{enumerate}
\item Suppose $b\in N_G(\csg{a})$ (the normalizer of $\csg{a}$ in $G$).  
        Then all out-neighbors of each element of $J$ are in $J$ and 
                  the induced subgraph of $\dpk{\exclude{\pi}}{G}$ on $J$ is $\dpk{\exclude{\pi}}{J}$.
         Note: If $\csg{a}$ is characteristic in $H$, such as when $H$ is cyclic, 
                  then $b\in N_G(\csg{a})$ for all $b\in K$.
\item \label{thm:complementedhall-fixed}
         Suppose $b\in C_G(\csg{a})$ (the centralizer of $\csg{a}$ in $G$).  
         If there is an edge from $b$ to $b'$ in $\dpk{}{K}$, 
         then for all $a'\in \gen{\csg{a}}$ there is an edge from $ab$ to $a' b'$ in $\dpk{\exclude{\pi}}{J}$.
         There are no other edges in $\dpk{\exclude{\pi}}{J}$.
         That is,  $\dqk{\exclude{\pi}}{J}\cong \dqk{}{K}$.
          Note: If  $G=H\times K$, then $b\in C_G(\csg{a})$ for all $b\in K$.
\end{enumerate}
\end{theorem}

\begin{proof} 
The multiplication in the inner semidirect product of $H$ by $K$ gives 
\[ (ab)^s = a\psi_{b}(a)\psi_{b}^2(a)\cdots \psi_{b}^{s-1}(a)b^s,\hbox{ where }\psi_b(a) = bab^{-1}.\]
    
 In (i), $\psi_b^k(a)\in\csg{a}$ for all $k$, 
       so $(ab)^s = a'' b^s$ for some $a''\in \csg{a}$ and $b^s\in \csg{b}$.
 Since $b\in N_G(\csg{a})$, $\csg{b}\subseteq N_G(\csg{a})$.
Thus, the out-neighbors of any vertex in $J$ is also in $J$.
The automorphism $\psi_b$ of $H$ (conjugation by $b$) maps characteristic subgroups to themselves, so such subgroups are normal in $G$.  Cyclic groups have a unique (characteristic) subgroup of each order.
 
In (ii), $\psi_b(a)  = a$, so $(ab)^s = a^sb^s$. Since $a$ and $b$ have coprime orders, the Chinese remainder theorem gives some $s$ such that $s$ is congruent to any specified value modulo $\oo{b}$ 
(we exclude values divisible by elements of $\pi$) and  any specified value modulo $\oo{a}$.  
Thus, there is some $s$ such that $b^s=b'$ and $a^s = a'$ as given in (ii), so $(ab)^s = a'b'$, as required.
When $G=H\times K$, elements of $H$ commute with the elements of $K$.
\end{proof} 

In the case $G=H\times K$, Theorem \ref{thm:complementedhall-fixed} gives an analog for $K\cong G/H$.

\begin{corollary}\label{cor:excludedindirectproduct}
Let  $H$ and $K$ be groups with coprime orders. Let $\pi$ be the set of prime divisors of $\ord{H}$.
Then $\dqk{\exclude{\pi}}{H\times K}$ is a disjoint union of copies of $\dqk{}{K}$, 
with one copy for each cyclic subgroup of $H$. (See Fig.~\ref{fig:Z12quotient}.)
\end{corollary}

\begin{example}
We illustrate Corollary \ref{cor:excludedindirectproduct} with  $\mathbb{Z}_{12}\cong\mathbb{Z}_3\times\mathbb{Z}_4$.
\begin{figure}[h]
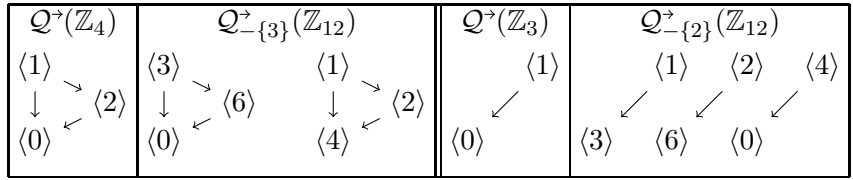

\[
\begin{array}{|c|c||c|c|}
\hline
\dqk{}{\mathbb{Z}_{4}} &\dqk{\exclude{\{3\}}}{\mathbb{Z}_{12}} &  \dqk{}{\mathbb{Z}_{3}} &  \dqk{\exclude{\{2\}}}{\mathbb{Z}_{12}} \\
\tikz
   \graph [no placement] {            
          1/$\csg{1}$ [x=-2,y=2] -> 2/$\csg{2}$ [x=-1, y=1.5];
          1->0/$\csg{0}$ [x=-2,y=1]; 
          2->0;
         };       
&     
\tikz
   \graph [no placement] {            
         1/$\csg{1}$ [x=.25,y=2] -> 2/$\csg{2}$ [x=1.25,y=1.5]; 
         2 -> 4/$\csg{4}$  [x=.25,y=1];
         1 -> 4;
          3/$\csg{3}$ [x=-2,y=2] -> 6/$\csg{6}$ [x=-1, y=1.5];
          3->0/$\csg{0}$ [x=-2,y=1]; 
          6->0;
         }; 
&
\tikz
   \graph [no placement] {            
         1/$\csg{1}$  [x=3,y=3]->0/$\csg{0}$ [x=2,y=2];
         };
&
\tikz
   \graph [no placement] {            
         1/$\csg{1}$ [x=1,y=3] -> 3/$\csg{3}$ [x=0,y=2];
         2/$\csg{2}$ [x=2,y=3]   -> 6/$\csg{6}$ [x=1, y=2];
         4/$\csg{4}$  [x=3,y=3]->0/$\csg{0}$ [x=2,y=2];
         };
         \\
\hline
\end{array}
\]
\caption{Example of Corollary \ref{cor:excludedindirectproduct} with 
      $\mathbb{Z}_{12}\cong\mathbb{Z}_3\times\mathbb{Z}_4$.}
\label{fig:Z12quotient}
\end{figure}
\end{example}

\begin{lemma}\label{lem:sdpZ2inv}
Let $H$ be an abelian group of finite order, and let $\pi$ be the set of prime divisors of $\ord{H}$. 
Let $K$ be a group with order coprime to that of $H$.  Let $\mathrm{id}$ and $\mathrm{inv}$ be the identity and inverse automorphisms of $H$, and let $\phi:K\rightarrow \{\mathrm{id}, \mathrm{inv}\}$ be a homomorphism.  
Consider $G=H\rtimes_\phi K$. Let $a\in H$, $b\in K$.  
\begin{enumerate}
\item  If $\phi_b=\mathrm{id}$, 
        then in $\dpk{\exclude{\pi}}{G}$, $(a,b)$ has out-neighbors 
        \[ (\gen{\csg{a}}\times \csg{b})\setminus\{(a,b)\}.\]  
 \item If $\phi_b=\mathrm{inv}$, 
        then in both $\dpk{\exclude{\pi}}{G}$ and  $\dpk{}{G}$,  $(a,b)$ has out-neighbors 
        \[\{(a,b^{2\ell+1}): 0< \ell<  \frac{\oo{b}}{2}\} \cup \{(0,b^{2\ell}):  0 \leq\ell<\frac{\oo{b}}{2}\}.\]  
 \end{enumerate}
 \end{lemma}

\begin{proof} 
If $\phi_b=\mathrm{id}$, then Theorem \ref{thm:complementedhall-fixed} (as applied to the outer construction of $G$) gives (i). 
If $\phi_b=\mathrm{inv}$, then compute $(a,b)^2=(a+\phi_b(a), b^2) = (a-a, b^2)=(0,b^2)$, $\ldots$, 
$(a,b)^{2\ell}=(0, b^{2\ell})$, $(a,b)^{2\ell+1}=(a, b^{2\ell+1})$.  Note that $b$ has even order since its homomorphic image does, so $|H|$ is odd.  Note that $(a,b)^{\oo{b}} = (0,0)$, so $\oo{(a,b)}=\oo{b}$. 
If the exponent $2\ell$ or $2\ell+1$ is divisible by elements of $\pi$,  then $2\ell+k\cdot\oo{b}$ or $2\ell+1+k\cdot\oo{b}$ is not for some value of $k$ while the parity remains unchanged.  Thus, none of these edges are excluded by $\pi$.
\end{proof} 

\begin{example}\label{ex:Z3xZ3spdZ2inv}
Let $G=(\mathbb{Z}_3\times\mathbb{Z}_3)\rtimes_\psi \mathbb{Z}_2$, where the homomorphism $\psi$ of $\mathbb{Z}_2$ into $\aut{\mathbb{Z}_3\times\mathbb{Z}_3}$ is given by $\psi_0$ is the identity and $\psi_1$ is the inverse map. 
($G$ is isomorphic to the generalized dihedral group for $\mathbb{Z}_3\times\mathbb{Z}_3$). 
Every element of the form $((a,b),1)$ is its own inverse.
Figure \ref{fig:Z3Z3Z2ex2} shows $\dpk{}{G}$, illustrating Lemma \ref{lem:sdpZ2inv}.
Thick arrows, representing powers  divisible by 3, are not present in $\dpk{\exclude{\{3\}}}{G}$. 

\begin{figure}[h]
\[
\tikz
   \graph [no placement] {            
          1a/$\scriptstyle{100}$ [x=-.8,y=1.5] <-> 2a/$\scriptstyle{200}$ [x=.8, y=1.5];
          5a/$\scriptstyle{110}$ [x=4.25,y=1.75] <->[bend left] 6a/$\scriptstyle{220}$ [x=4.25, y=-1.75];
          7a/$\scriptstyle{120}$ [x=-1.25,y=.8] <-> 8a/$\scriptstyle{210}$ [x=-1.25, y=-.8];
           3a/$\scriptstyle{010}$ [x=-.8,y=-1.5] <-> 4a/$\scriptstyle{020}$ [x=.8, y=-1.5];
           {1a, 2a,3a,4a,7a,8a}->[ultra thick] 0a/$\scriptstyle{000}$;
           5a ->[bend right, ultra thick] 0a;  6a ->[bend left, ultra thick] 0a;
            5c/$\scriptstyle{111}$ [x=2.5,y=.4] ;  6c/$\scriptstyle{221}$ [x=2.5, y=-.4];
            0c/$\scriptstyle{001}$ [x=1.5, y=0] -> 0a;
           1c/$\scriptstyle{101}$ [x=3,y=1.5];  2c/$\scriptstyle{021}$ [x=3, y=.8];
           7c/$\scriptstyle{211}$ [x=1.5,y=.8]->0a; 
           3c/$\scriptstyle{011}$ [x=3,y=-1.5]; 4c/$\scriptstyle{201}$ [x=3, y=-.8]; 
           8c/$\scriptstyle{121}$ [x=1.5,y=-.8]->0a;  
            1c ->[bend right] 0a;
            3c ->[bend left] 0a;
            2c-> 0a;
            4c-> 0a;
            {5c,6c}-> 0a;
         };
\tikz
   \graph [no placement] {            
          l/$\dqk{\exclude{\{3\}}}{(\mathbb{Z}_3\times\mathbb{Z}_3)\rtimes_\psi \mathbb{Z}_2}$ [x=1.25, y=2];
          1a/$\scriptstyle{\csg{100}}$ [x=0,y=1] ;
          5a/$\scriptstyle{\csg{110}}$ [x=2.5,y=0];
          7a/$\scriptstyle{\csg{120}}$ [x=-.75,y=0] ;
          3a/$\scriptstyle{\csg{010}}$ [x=00,y=-1] ;        
          0a/$\scriptstyle{\csg{000}}$;
          0c/$\scriptstyle{\csg{001}}$ [x=1.25, y=0] -> 0a;
          7c/$\scriptstyle{\csg{211}}$ [x=1.25,y=1] -> 0a; 
          8c/$\scriptstyle{\csg{121}}$ [x=1.25,y=-1] ->0a;     
          5c/$\scriptstyle{\csg{111}}$ [x=2,y=.5] ;
          1c/$\scriptstyle{\csg{101}}$ [x=2,y=1];
          3c/$\scriptstyle{\csg{011}}$ [x=2,y=-1];  
          5c-> 0a; 
          1c->0a; 
          3c-> 0a;  
           };
 \]
 \caption{The power graph of a semidirect product $(\mathbb{Z}_3\times\mathbb{Z}_3)\rtimes_\psi \mathbb{Z}_2$.}
 \label{fig:Z3Z3Z2ex2}
 \end{figure}          
\end{example}

Compare Example \ref{ex:Z3xZ3spdZ2inv} to another semidirect product of the same groups.

\begin{example}
Let $G=(\mathbb{Z}_3\times\mathbb{Z}_3)\rtimes_\phi \mathbb{Z}_2$, where the homomorphism $\phi$ of $\mathbb{Z}_2$ into $\aut{\mathbb{Z}_3\times\mathbb{Z}_3}$ is given by $\phi_0$ is the identity and $\phi_1$ swaps components. 
 ($G$ is isomorphic to $S_3\times \mathbb{Z}_3$).
Figure \ref{fig:Z3Z3Z2ex1} shows $\dpk{}{G}$.
Thick arrows, representing powers  divisible by 3, are not present in $\dpk{\exclude{\{3\}}}{G}$. 
Elements joined by a bidirectional edge generate the same cyclic subgroup and represent the same element of the quotient power graph.  

\begin{figure}[h]
\[
\tikz
   \graph [no placement] {            
          1a/$\scriptstyle{100}$ [x=-.8,y=1.5] <-> 2a/$\scriptstyle{200}$ [x=.8, y=1.5];
          5a/$\scriptstyle{110}$ [x=4.25,y=1.75] <->[bend left] 6a/$\scriptstyle{220}$ [x=4.25, y=-1.75];
          7a/$\scriptstyle{120}$ [x=-1.25,y=.8] <-> 8a/$\scriptstyle{210}$ [x=-1.25, y=-.8];
           3a/$\scriptstyle{010}$ [x=-.8,y=-1.5] <-> 4a/$\scriptstyle{020}$ [x=.8, y=-1.5];
           {1a, 2a,3a,4a,7a,8a}->[ultra thick] 0a/$\scriptstyle{000}$;
           5a ->[bend right, ultra thick] 0a;  6a ->[bend left, ultra thick] 0a;
             5c/$\scriptstyle{111}$ [x=2.5,y=.4] <-> 6c/$\scriptstyle{221}$ [x=2.5, y=-.4];
             0c/$\scriptstyle{001}$ [x=1.5, y=0] -> 0a;
             5c->[bend right] 5a; 6c->[bend right] 5a;
             5c->[bend left] 6a; 6c->[bend left] 6a;
             5c->[ultra thick] 0c; 6c->[ultra thick] 0c; 
           1c/$\scriptstyle{101}$ [x=3,y=1.5] <-> 2c/$\scriptstyle{021}$ [x=3, y=.8];
           7c/$\scriptstyle{211}$ [x=1.5,y=.8]; 
           3c/$\scriptstyle{011}$ [x=3,y=-1.5] <-> 4c/$\scriptstyle{201}$ [x=3, y=-.8]; 
           8c/$\scriptstyle{121}$ [x=1.5,y=-.8];  
            1c->5a; 1c->[bend left] 6a; 2c->[bend right] 5a; 2c->[bend left] 6a;  
            3c->[bend right] 5a; 3c->6a;  4c->[bend right] 5a; 4c->[bend left] 6a;
            1c->[ultra thick] 7c; 2c->[ultra thick]7c; 7c->0a;
            3c->[ultra thick] 8c; 4c->[ultra thick] 8c; 8c->0a;
            1c ->[bend right, ultra thick] 0a;
            3c ->[bend left, ultra thick] 0a;
            2c->[ultra thick] 0a;
            4c->[ultra thick] 0a;
            {5c,6c}->[ultra thick] 0a;
         };
\tikz
   \graph [no placement] {            
          l/$\dqk{\exclude{\{3\}}}{(\mathbb{Z}_3\times\mathbb{Z}_3)\rtimes_\phi \mathbb{Z}_2}$ [x=1.25, y=2];
          1a/$\scriptstyle{\csg{100}}$ [x=0,y=1] ;
          5a/$\scriptstyle{\csg{110}}$ [x=3.25,y=0];
          7a/$\scriptstyle{\csg{120}}$ [x=-.75,y=0] ;
          3a/$\scriptstyle{\csg{010}}$ [x=00,y=-1] ;        
          0a/$\scriptstyle{\csg{000}}$;
          0c/$\scriptstyle{\csg{001}}$ [x=1.25, y=0] -> 0a;
          7c/$\scriptstyle{\csg{211}}$ [x=1.25,y=1] -> 0a; 
          8c/$\scriptstyle{\csg{121}}$ [x=1.25,y=-1] ->0a;     
          5c/$\scriptstyle{\csg{111}}$ [x=2,y=0] ;
          1c/$\scriptstyle{\csg{101}}$ [x=2,y=1];
          3c/$\scriptstyle{\csg{011}}$ [x=2,y=-1];  
          5c-> 5a; 
          1c->5a; 
          3c-> 5a;  
           };
 \]
 \caption{The power graph of another semidirect product $(\mathbb{Z}_3\times\mathbb{Z}_3)\rtimes_\psi \mathbb{Z}_2$.}
 \label{fig:Z3Z3Z2ex1}
 \end{figure}              
\end{example}

Power graphs of nilpotent groups have been studied in a number of papers, including \cite{MR4469917,MR4504122,MR4207766,MR4746158,MR4835924,MR4718723,MR4098926}.
The power graph of a nilpotent group with a normal Hall subgroup is discussed in \cite{MR4391495}.
Power graphs of abelian groups and $p$-groups have been studied as well, for example in  \cite{MR4977007,MR5020628,MR4223622,MR4389908}.
Here we see that excluded power graphs of nilpotent groups are nice. 

\begin{corollary}\label{cor:nilpotentexclucded}
Let $G$ be a finite nilpotent group. Let $\pi$ be the set of prime divisors of $\ord{G}$, and for each $p\in \pi$, let $S_p$ be the Sylow $p$-subgroup of $G$.  
Let $\rho$ be any proper nonempty subset of $\pi$, and let $\bar{\rho} = \pi\setminus \rho$.
Then ${\dqk{\exclude{\rho}}{G}\cong \dqk{}{\prod_{p\in \bar{\rho}} S_p}}$.
\end{corollary}

\begin{proof} 
A finite nilpotent group is the direct product of its Sylow subgroups.
Induct on the size of $\pi$.  
Theorem \ref{thm:complementedhall-fixed} treats the case where $\rho$ has size 1.
Suppose ${\dqk{\exclude{\rho}}{G}\cong \dqk{}{\prod_{p\in \bar{\rho}} S_p}}$ for any proper nonempty subset $\rho$ of $\pi$.  Suppose there is a prime $q\in\pi$ such that $\rho'=\rho\cup\{q\}$ is proper.
Then applying Theorem \ref{thm:complementedhall-fixed} gives
${\dqk{\exclude{\rho'}}{G}\cong \dqk{\exclude{\{q\}}}\prod_{p\in \bar{\rho}} S_p}\cong \dqk{}{\prod_{p\in \bar{\rho}'} S_p}$.
\end{proof} 

While solvable groups have Hall subgroups of all possible orders, the Hall subgroups are generally not complemented (see \cite{MR1307623}, for example). 
Thus, the results of this section, which exploit the semidirect product construction, are not directly applicable in general solvable groups.  
The Fitting subgroup, being nilpotent, will contribute a substructure described in Corollary \ref{cor:nilpotentexclucded}.
Nonetheless, we hope that the structure of an excluded power graph of a solvable group  will reflect the presence of Hall subgroups in a meaningful way. 
The results concerning power graphs of solvable groups, such as \cite{MR3691533,MR4335775,MR4764545,MR4835924}, are less extensive than those for nilpotent groups.

\section{Directed Cliques}\label{sec:disjointcliques}

A finite group has complete undirected power graph if and only if the group is cyclic with prime power order \cite{MR2511776}. 
We give an analog for $\pi$-excluded power graphs for sets of primes $\pi$, which turns out to involve semidirect products.  
Lemma \ref{lem:excludeallprimedivisors} describes digraphs which are the disjoint union of directed cliques.

\begin{lemma}\label{lem:undirunionofcliques}
Let $\pi$ be a set of primes, and let $q$ be a prime not in $\pi$. 
Let $P$ be a finite $\pi$-group, and let $s$ be a nonnegative integer.
Fix a homomorphism $\phi:P\rightarrow\aut{ \mathbb{Z}_{q^s}}$.
Then $\pk{\exclude{\pi}}{\mathbb{Z}_{q^s}\rtimes_\phi P}$ is a disjoint union of cliques.
\end{lemma}

\begin{proof} 
Abbreviate $G=\mathbb{Z}_{q^s}\rtimes_\phi P$.  Say $P$ has order $n$, where the prime factors of $n$ lie in $\pi$.
Write the elements of $G$ as ordered pairs $(a,b)$ with $a\in \mathbb{Z}_{q^s}$, $b\in P$. 
Use additive notation for $\mathbb{Z}_{q^s}$, writing $n\cdot a$ for the sum of $n$ many $a$'s.
Each subgroup of $\mathbb{Z}_{q^s}$ consists of the multiples of some generator of the subgroup.
Use multiplicative notation for $P$ and $G$, writing $b^n$ for the product of $n$ many $b$'s.
Write $\phi_b\in \aut{ \mathbb{Z}_{q^s}}$ to be the image of $b$ under $\phi$.

We first consider edges in $\dpk{\exclude{\pi}}{G}$.
Note that $(a,b)^k = (a + \phi_{b}(a) + \cdots + \phi_{b^{k-1}}(a), b^k)$.
Exponents $k$ divisible any elements of $\pi$ are excluded, so $b^{k}\in\gen{\csg{b}}$ 
for all powers $k$ used to construct the edges of $\dpk{\exclude{\pi}}{G}$.

We claim that if $\phi_b$ is the identity automorphism, 
then the Cartesian product $\mathbb{Z}_{q^s}\times \gen{\csg{b}}$ is not only a clique in $\pk{\exclude{\pi}}{G}$, but also a connected component.
Here $(a,b)^k =(k \cdot a, b^k)$.   
No edges enter or leave the set in question since $b^{k}\in\gen{\csg{b}}$.
The order of $a$ is a power of $q$ and the order of $b$ has all factors in $\pi$, so these orders are coprime.  
By the Chinese remainder theorem, there is a solution to $k\equiv u \pmod {\oo{a}}$ and $k\equiv v \pmod {\oo{b}}$, 
so there is an edge in $\dpk{\exclude{\pi}}{G}$ from $(a,b)$ to $(a',b')$ for all $a'\in\csg{a}$ and $b'\in\gen{\csg{b}}$.  
The subgroups of a cyclic group  of prime power order (such as $\csg{a}$) are linearly ordered, 
so for any pair of elements, at least one is a power of the other.  The claim follows.  
This claim gives the result in the case that $G=\mathbb{Z}_{q^s}\times P$.  
(This case also follows from the fact that $\pk{}{\mathbb{Z}_{q^s}}$ is a clique \cite{MR2511776} and Corollary \ref{cor:excludedindirectproduct}).   

Now consider general $\phi_b$.  Here $b\mapsto \phi_b$ is a homomorphism from $P$ to $\aut{\mathbb{Z}_{q^s}}$, and these groups have respective orders $n=p_1^{t_1}\ldots p_h^{t_h}$ and $\varphi(q^s)=(q-1)q^{s-1}$.  
For $b\in P$, let $m_b=\oo{\phi_b}$.  By Lagrange's Theorem,  
        $m_b \mid \oo{b}$, so   $m_b=p_1^{r_1} \cdots p_h^{r_h}$ for some $e_i$ $(0\leq r_i\leq t_\ell, 1\leq i \leq h)$.
 Additionally,        
        $m_b \mid \ord{\aut{\mathbb{Z}_{q^s}}}=\varphi(q^s)= (q-1)q^{s-1}$.  
In particular, $p_i^{r_i}\mid q-1$.  

Every automorphism of $\mathbb{Z}_{q^t}$ group is determined by the image of $1$.
The case $\phi_{b}(1) = 1$ ($\phi_b$ is the identity automorphism) is treated above.
Suppose $\phi_{b}(1) = c_b$ for some $c_b\in\gen{\mathbb{Z}_{q^s}}\setminus\{1\}$.  
Then $\phi_b(a) =  c_b \cdot a$.   
By the homomorphism property, $\phi_{b^k}(a)=\phi_{b}^k(a) = c_b^k \cdot a$. 
Now $a + \phi_{b}(a) + \cdots + \phi_{b^{k-1}}(a) = a +c_b\cdot a+ c_b^{2} \cdot a + \cdots + c_b^{k-1}\cdot a
                  = (c_b^k-1)/(c_b-1)\cdot a $.
The multiplicative order of $c_b$ mod $q^s$ is $\oo{\phi_b}=m_b$. 
Thus, $(c_b^k-1)/(c_b-1)$ takes on $m_b$-many distinct values modulo $q^s$.
Recall that we are excluding $k$ which are multiples of $p$, and 
we do not consider $k=1$ (as there are no loops). 

We have edges in $\dpk{\exclude{\pi}}{G}$
       from $(a,b)$ to $((c_b^k-1)/(c_b-1)\cdot a, b^k)$ for $1 < k < m_b$ with $k$ not divisible by any element of $\pi$.
For such a $k$,  consider the edges from $(a(c_b^k-1)/(c_b-1), b^k)$.
Note that $\phi_{b^k} = \phi_b^k$, so $\phi_{b^k}(a) = c_b^k a$.
Taking $\ell\not \in \pi$, we compute using the above:
 \[
 (\frac{c_b^k-1}{c_b-1}\cdot a, b^k)^\ell 
     = (\frac{c_b^k-1}{c_b-1} \cdot \frac{c_{b^k}^\ell-1}{c_{b^k}-1}\cdot a , {b^{k}}^\ell)
     =  (\frac{c_{b}^{k\ell}-1}{c_b-1}\cdot a , b^{k\ell}) = (a,b)^{k\ell}.
\]
Note that $\ell$ is also allowed to take any value such that $1\leq \ell < m_b$ with $\ell$ not divisible by any element of $\pi$.
Thus, 
     $\{ (a,b)\}\cup  \{ (a(c_b^k-1)/(c_b-1), b^k) : 1 < k < m_b,\,\hbox{$k$ not divisible by any element of $\pi$}\}$ 
is a directed clique in $\dpk{\exclude{\pi}}{G}$. 

Suppose there is an edge from both $(a', b')$ and $(a'', b'')$ to the same element $(a, b)$.  
Then by the above there is an edge from $(a, b)$ to both $(a', b')$ and $(a'', b'')$, so there  are edges in each direction between $(a', b')$ and $(a'', b'')$.  
That is,  $(a', b')$ and $(a'', b'')$ are in the directed clique described above.
\end{proof} 

\begin{corollary}\label{cor:dpk-sdp-ZP}
Let $\pi$ be a set of primes, and let $q$ be a prime not in $\pi$. 
Let $P$ be a $\pi$-group, and let $s$ be a nonnegative integer.
Fix a homomorphism $\phi:P\rightarrow\aut{ \mathbb{Z}_{q^s}}$, and 
let $G=\mathbb{Z}_{q^s}\rtimes_\phi P$.
Pick $b\in P$, and let $c=\phi_b(1)$.
\begin{enumerate}
\item If $c=1$, 
         then $\mathbb{Z}_{q^s}\times \gen{\csg{b}}$ is a weakly connected component of $\dpk{\exclude{\pi}}{G}$.
         If there is an edge from $a$ to $a'$ in $\dpk{}{\mathbb{Z}_{q^s}}$, 
         then there is an edge from $(a, b)$ to $(a', b')$ for all $b'\in\gen{\csg{b}}$.
\item If $c\not=1$,  then for all $a\in \mathbb{Z}_{q^s}$, 
                $\{ (a,b)\}\cup  \{ (a(c^k-1)/(c-1), b^k) : 1 < k < m_b,\,\hbox{$k$ not divisible by any element of $\pi$}\}$ 
         is a weakly connected component of and a directed clique in $\dpk{\exclude{\pi}}{G}$.
\end{enumerate}
\end{corollary}

\begin{example}
Consider
$\dpk{\exclude{\{3\}}}{\mathbb{Z}_7\rtimes_\phi\mathbb{Z}_3}$, where $\phi_b(a)=a\cdot 2^b$ (Fig.~\ref{fig:dpg(Z7sdpZ3)}). Here $c_1=2$ has multiplicative order 3. The only multiplier of interest is
 $(2^2-1)/(2-1)=3$, so each element $(a,1)$ belongs to a disjoint clique with $(3a,2)$, where $3a$ is taken mod 7. This illustrates Lemma \ref{lem:undirunionofcliques} and Corollary \ref{cor:dpk-sdp-ZP}.

\begin{figure}[h]
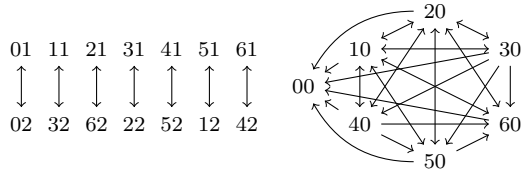

\[
\tikz
   \graph [no placement] {
          01/$\scriptstyle{01}$ [x=-2.5,y=.5]  <-> 02/$\scriptstyle{02}$  [x=-2.5,y=-.5];            
          11/$\scriptstyle{11}$ [x=-2,y=.5]  <-> 32/$\scriptstyle{32}$  [x=-2,y=-.5];
          21/$\scriptstyle{21}$ [x=-1.5,y=.5]  <-> 62/$\scriptstyle{62}$  [x=-1.5,y=-.5];
          31/$\scriptstyle{31}$ [x=-1,y=.5]  <-> 22/$\scriptstyle{22}$  [x=-1,y=-.5];
          41/$\scriptstyle{41}$ [x=-.5,y=.5]  <-> 52/$\scriptstyle{52}$  [x=-.5,y=-.5];
          51/$\scriptstyle{51}$ [x=0,y=.5]  <-> 12/$\scriptstyle{12}$  [x=0,y=-.5];
          61/$\scriptstyle{61}$ [x=.5,y=.5]  <-> 42/$\scriptstyle{42}$  [x=.5,y=-.5];
          00/$\scriptstyle{00}$ [x=1.25,y=0] ;
          10/$\scriptstyle{10}$ [x=2,y=.5] -> 00;
          20/$\scriptstyle{20}$ [x=3,y=1] ->[bend right] 00;
          30/$\scriptstyle{30}$ [x=4,y=.5] -> 00; 
          40/$\scriptstyle{40}$ [x=2,y=-.5] -> 00;
          50/$\scriptstyle{50}$ [x=3,y=-1] ->[bend left] 00;
          60/$\scriptstyle{60}$ [x=4,y=-.5] -> 00; 
          10<->{20,30,40,50,60};
          20<->{30,40,50,60};
          30->{40,50,60};
          40->{50,60};
          50->{60};
           };
\]
\caption{The $\{3\}$-excluded power graph of a semidirect product $\mathbb{Z}_7\rtimes_\phi\mathbb{Z}_3$.}
\label{fig:dpg(Z7sdpZ3)}
\end{figure}
\end{example}

\begin{theorem}\label{thm:unioncomplete}
Let $G$ be a group, and let $H$ be a subgroup of $G$.   Let $\pi$ be a set of primes.     
\begin{enumerate}
\item \label{thm:unioncomplete-dir}
The following are equivalent.
\begin{enumerate}
\item  $\dpk{\exclude{\pi}}{H}$ is a disjoint union of directed cliques.
\item All prime divisors of $|H|$ are in $\pi$.
\end{enumerate}
\item \label{thm:unioncomplete-undir} 
The following are equivalent.
\begin{enumerate}
\item $\pk{\exclude{\pi}}{H}$ is a disjoint union of cliques. 
\item $H\cong \mathbb{Z}_{q^s}\rtimes P$ for some prime $q\not\in\pi$, nonnegative integer $s$, and $\pi$-group $P$.
\end{enumerate}
\end{enumerate}
\end{theorem}

\begin{proof} 
(i): Suppose  $\ord{H}$ is divisible by some prime $q$ not in  $\pi$.  
Then there is an edge from each non-identity element of $q$-power order to the identity, but no edge in the reverse direction.    
Thus, $\dpk{\exclude{\pi}}{H}$ is not a disjoint union of directed cliques. Hence, (ia) implies (ib).  
Conversely, note that $\dpk{\exclude{\pi}}{P}$ is a disjoint union of directed cliques by Lemma \ref{lem:excludeallprimedivisors}, so (ib) implies (ia).

(ii): First suppose $\pk{\exclude{\pi}}{H}$ is a disjoint union of cliques.
Suppose $\ord{H}$ is divisible by distinct primes $q$ and $r$ not in $\pi$.  
Then there are cyclic subgroups $Q$ of order $q$ and $R$ of order $r$.  
There is an edge from each element of $Q$ and $R$ to $\id$ 
           but no edges between non-identity elements of $Q$ and $R$.  
Hence, $\pk{\exclude{\pi}}{H}$ is not the union of (directed) cliques in this case. 
Thus, $|H|$ has at most one prime factor not in $\pi$.

If $H$ is a $\pi$-group, then  (iib) holds with $s=0$.  
Suppose $\ord{H}$ is divisible by one prime $q$ not in $\pi$ and  that the subgroups with $q$-power order are not linearly ordered. 
Then there are elements $x$, $y\in H$ with $q$-power order, neither of which is a power of the other, while the identity is a $q$-power of both.  
So $\pk{\exclude{\pi}}{H}$ is not the disjoint union of cliques unless the elements of $q$-power order form a cyclic subgroup $\mathbb{Z}_{q^s}$ (the case where the subgroups are linearly ordered).  
This $\mathbb{Z}_{q^s}$ is characteristic, so it is a normal Sylow subgroup, and hence is complemented by the Schur-Zassenhaus theorem.  
Therefore, $G\cong \mathbb{Z}_{q^s}\rtimes P$ for some $\pi$-group $P$.  Thus, (iia) implies (iib).
Conversely, note that $\pk{\exclude{\pi}}{\mathbb{Z}_{q^s}\rtimes P}$ is a disjoint union of directed cliques by Lemma \ref{lem:undirunionofcliques}, so (iib) implies (iia).
\end{proof} 

\begin{example}
The dihedral group $D_{2n}$ is isomorphic to $\mathbb{Z}_n\rtimes \mathbb{Z}_{2}$.   
Theorem \ref{thm:unioncomplete-undir} applies to the case $n=q^a$ for some odd prime $q$. (If $q=2$, then $D_{2n}$ is a $2$-group and subject to Theorem \ref{thm:unioncomplete-dir}).
In  $\pk{\exclude{\{2\}}}{D_{2q^{a}}}$, the reflections are isolated vertices and the subgroup of rotations forms a single connected component. 
As a non-example, in $\pk{\exclude{\{q\}}}{D_{2q^{a}}}$, the reflections form a star centered around the identity, and the non-identity rotations form a clique.
\end{example}


\end{document}